\documentclass[12pt]{article}

\usepackage{amsmath, amsthm, amssymb, amsfonts, enumerate, dsfont}
\usepackage[applemac]{inputenc} 

\def\dis{\displaystyle}
\def \R{I\!\!R}

\newtheorem{thm}{Theorem}[section]
\newtheorem{cor}[thm]{Corollary}
\newtheorem{lem}[thm]{Lemma}
\newtheorem{pro}[thm]{Proposition}
\newtheorem{defi}[thm]{Definition}

\newtheorem{exa}[thm]{Example}

\def \lim{ {\rm lim} }

\def \1{\mathbf{1}}

\def\abstract{\begin{center} \small\bf Abstract\end{center}\small}

 \title{General limit value in Dynamic Programming}

 \author{J\'er\^ome Renault\thanks{%
TSE (GREMAQ, Universit\' e Toulouse 1 Capitole), 21 all\' ee de Brienne, 31000
Toulouse, France. E-mail: \textsf{jerome.renault@tse-fr.eu}.   The author  gratefully acknowledges  the support of the Agence Nationale de la Recherche, under grant ANR JEUDY, ANR-10-BLAN 0112.
  }  \date{\today}}

\begin{document}

\maketitle


\begin{abstract}  We consider a   dynamic programming problem  with arbitrary state space and bounded rewards.  Is it possible to define in an unique way a   limit value for the problem, where the ``patience" of the decision-maker tends to infinity ? We  consider, for  each evaluation $\theta$ (a  probability  distribution  over  positive integers)  the value function $v_{\theta}$ of the problem where the weight of any stage $t$ is given by $\theta_t$, and we investigate the uniform convergence of  a sequence  $(v_{\theta^k})_k$ when  the ``impatience" of the evaluations vanishes, in the sense that $\sum_{t} | \theta^k_{t}-\theta^k_{t+1}| \rightarrow_{k \to \infty} 0$.   We prove that this uniform convergence happens if and only if the metric space $\{v_{\theta^k}, k\geq 1\}$ is totally bounded. Moreover  there exists a particular function $v^*$, independent of the particular chosen  sequence $({\theta^k})_k$,  such that    any limit point of such sequence of value functions is  precisely $v^*$.  Consequently, while speaking of uniform convergence of the value functions, $v^*$ may be considered as the unique possible  limit  when the patience of the decision-maker tends to infinity. The result applies in particular  to discounted payoffs when the discount factor vanishes, as well as to average payoffs where the number of stages goes to infinity, and also to models  with  stochastic transitions. We present  tractable corollaries, and   we discuss counterexamples and a conjecture. \\

 \end{abstract}

\section{Introduction}

We consider a   dynamic programming problem  with arbitrary state space $Z$ and bounded rewards.  Is it possible to define in an unique way a possible  limit value for the problem, where the ``patience" of the decision-maker tends to infinity ? 

For  each evaluation (probability  distribution  over  positive integers) $\theta=(\theta_t)_{t\geq 1}$, we consider   the value function $v_{\theta}$ of the problem where the  initial state is arbitrary in $Z$ and  the weight of any stage $t$ is given by $\theta_t$. The  {\it total variation} of $\theta$, that we also wall the {\it impatience} of $\theta$, is defined by   $TV(\theta)=\sum_{t=1}^\infty |\theta_{t+1}-\theta_t|.$ For instance, for each positive integer $n$  the evaluation $\theta=(1/n,...,1/n,0,...,0,...)$ induces the value function $\bar{v}_n$ corresponding to the maximization of the mean  payoff for  the first $n$ stages; and for any $\lambda$ in $(0,1]$ the evaluation $\theta=(\lambda(1-\lambda)^{t-1})_t$ induces the discounted value function $v_\lambda$. 

A well known theorem of Hardy and Littlewood (see e.g. Lippman, 1969) implies that for an {\it uncontrolled} problem, the pointwise convergence of $(\bar{v}_n)_n$, when $n$ goes to infinity,  and of $(\bar{v}_\lambda)_\lambda$, when $\lambda$ goes to 0, are equivalent, and that in case of convergence both limits are the same. However, Lehrer and Sorin (1992) provided an example of a dynamic programming problem where $(\bar{v}_n)_n$ and $(\bar{v}_\lambda)_\lambda$ have different pointwise limits. But they also proved that the {\it uniform} convergence of $(\bar{v}_n)_n$  and of $(\bar{v}_\lambda)_\lambda$  are equivalent, with equality of the limit in case of convergence.  And Sorin and Monderer (1993) extended this result to   families of evaluations satisfying some conditions.  Mertens and Neyman (1982) proved that when the family $(\bar{v}_\lambda)_\lambda$ not only uniformly converges but has bounded variation, then the dynamic programming problem has a {\it uniform value}, in the sense  that for all initial state $z$ and $\varepsilon>0$, there exists a play with mean payoffs from stage 1 to stage $T$  at least $v-\varepsilon$ provided $T$ is large enough  (see also Lehrer Monderer 1994 and Sorin Moderer 1993 for   proofs that the uniform convergence of $(\bar{v}_\lambda)_\lambda$ or $(\bar{v}_n)_n$ does not imply the existence of the {\it uniform value} of the problem).   In this case of existence of a uniform value, one can show that all value functions $v_{\theta}$ are close to $v^*$, whenever  $\theta$ is a {\it non increasing} evaluation with small $\theta_1$. The reason is that whenever $\theta$ is non increasing, the $\theta$-payoff of a play can be expressed as a convex combination of the Cesàro values $(\bar{v}_n)_n$.  \\

In the present paper,  we investigate the uniform convergence of   sequences  $(v_{\theta^k})_k$ when  the ``impatience" of the evaluations vanishes, in the sense that $\sum_{t} | \theta^k_{t}-\theta^k_{t+1}| \rightarrow_{k \to \infty} 0$.   We will prove in theorem \ref{thm1} that this uniform convergence happens if and only if the metric space $\{v_{\theta^k}, k\geq 1\}$ (with the distance between functions given by the sup of their differences),  is totally bounded. Moreover  the uniform limit, whenever it exists, can only be the following   function, which is independent of the particular chosen  sequence $({\theta^k})_k$:
$$v^*=\inf_{\theta \in \Theta} \; \sup_{m \geq 0} \;  {{v}}_{m, \theta},$$ where  for each evaluation $\theta=(\theta_t)_{t\geq 1}$, the evaluation $m, \theta$ is defined as the evaluation with weight 0 for the first 
$m$ stages and with weight $\theta_{t-m}$ for stages $t>m$.  Consequently, while speaking of uniform convergence of the value functions when the patience of the decision-maker tends to infinity, $v^*$ can  be considered as the unique possible  limit value.   We also  give simple conditions on the state space, the payoffs and the transitions (mainly compactness, continuity and non expansiveness) implying the uniform convergence of such value functions. \\

The paper is organized as follows: section 2 contains the model and the main results, which  are shown to extend to the case of stochastic transitions. Section 3 contains a few  examples and counterexamples and  section 4 contains the proof of theorem \ref{thm1}. In the last section we formulate the following conjecture, which is shown to be true for uncontrolled problems:  does the uniform convergence of  $(\bar{v}_n)_n$, or equivalently of $(\bar{v}_\lambda)_\lambda$, implies the general convergence of the value functions, in the sense that: $\forall \varepsilon >0, \exists \alpha >0, \forall \theta \in \Theta \; s.t. \; TV(\theta)\leq \alpha, \; \|v_\theta-v^*\| \leq \varepsilon$ ?

\section{Model and results}

\subsection{General values in dynamic  programming problems}We consider a dynamic programming problem given by   a non empty set of states $Z$,   a correspondence   $F$ with non empty values from $Z$ to $Z$, and a mapping $r$ from $Z$ to $[0,1]$. 
   $Z$ is called the set of states, $F$ is the transition correspondence and $r$ is the reward (or payoff) function. An initial state $z_0$ in $Z$ defines    the following  dynamic programming problem:   a decision maker, also called player,  first has to select a new state $z_1$ in $F(z_0)$, and is rewarded by $r(z_1)$. Then he has to choose $z_2$ in $F(z_1)$, has a payoff of $r(z_2)$, etc... The  decision maker is interested in maximizing  his ``long-term" payoffs, for whatever it means.   From now on we fix $\Gamma=(Z,F,r)$, and for every state $z_0$ we denote by $\Gamma(z_0)=(Z,F,r,z_0)$ the corresponding problem  with initial state $z_0$. 
For $z_0$ in $Z$, a play  at $z_0$ is a sequence $s=(z_1,...,z_t,...) \in Z^{\infty}$ such that: $\forall t\geq 1, z_t \in F(z_{t-1})$. We denote by $S(z_0)$ the set of plays  at $z_0$, and by $S=\cup_{z_0\in Z} S(z_0)$ the set of all plays.  The set of bounded  functions from $Z$ to $\R$ is denoted by ${\cal V}$, and for $v$ and $v'$ in ${\cal V}$ we use the distance $d_\infty(v,v')= \sup_{z \in Z} |v(z)-v'(z)|$. \\

\noindent{\bf \it Ces\`{a}ro values.} For $n\geq 1$ and $s=(z_t)_{t\geq 1}\in S$,   the average payoff of the play $s$ up to stage $n$ is defined by:
${\gamma}_{\overline{n}}(s)= \frac{1}{n} \sum_{t=1}^n r(z_t).$
 And     the $n$-stage average  value of   $\Gamma(z_0)$ is:
$\dis  {v}_{\overline{n}}(z_0)= \sup_{s \in S(z_0)} \gamma_{\overline{n}}(s).$ By the   Bellman-Shapley  recursive formula,    for all $n$ and $z$  we have:  $n \;   {v}_{{\overline{n}}} (z)= \sup_{z'\in F(z)} \left( r(z')+ (n-1) \; {v}_{\overline{n-1}} (z')\right).$  We also  have  $|{v}_{{\overline{n}}}(z) -\sup_{z' \in F(z)}  {v}_{{\overline{n}}} (z')|\leq \frac{2}{n}$, and a pointwise limit of $( {v}_{\overline{n}})_n$ should satisfy $ v(z)= \sup_{z' \in F(z)} v (z')$ for all $z$. \\

\noindent{\bf \it Discounted  values.} Given $\lambda\in (0,1]$, the $\lambda$-discounted payoff of a play $s=(z_t)_{t\geq 1}$ is $\gamma_\lambda(s)= \lambda \sum_{t=1}^\infty (1-\lambda)^{t-1} r(z_t)$, and the  $\lambda$-discounted value at the initial state $z_0$ is $v_\lambda(z_0)= \sup_{s \in S(z_0)} \gamma_\lambda(s).$ It is easily proved that $v_\lambda$ is the unique mapping in ${\cal V}$  satisyfing the fixed point equation : $\forall z\in Z,   v_{\lambda} (z)= \sup_{z'\in F(z)} \left( \lambda \; r(z')+ (1-\lambda) \; v_{\lambda}(z')\right).$ It implies $|v_{\lambda}(z) -\sup_{z' \in F(z)} v_{\lambda} (z')|\leq \lambda$, and a  pointwise limit of $(v_\lambda)_\lambda$ should also satisfy $v(z)= \sup_{z' \in F(z)} v (z')$ for all $z$. \\

\noindent{\bf \it General values.} We denote by $\Theta$ the set of probability distributions  over positive integers.
An element  $\theta=(\theta_t)_{t \geq 1}$ in $\Theta$ is called an evaluation.

\begin{defi}    $\;$ 

The $\theta$-payoff of a play $s=(z_t)_{t\geq 1}$ is  
$\gamma_\theta(s)=   \sum_{t=1}^\infty \theta_tr(z_t),$

\vspace{0,3cm}

  and   the $\theta$-value of   $\Gamma(z_0)$ is 
$\dis v_\theta(z_0)= \sup_{s \in S(z_0)} \gamma_\theta(s).$ \end{defi}

\indent For each stage $t$ we denote by $\delta_t$ the Dirac mass on stage $t$ and by $\overline{n}$   the Ces\`{a}ro evaluation $(1/n,...,1/n,0,...,0,...)=(1/n) \sum_{t=1}^n \delta_t$, so that the notation $v_\theta$ for $\theta=\overline{n}$ coincide with the Ces\`{a}ro-value $v_{\overline n}$, also written $\bar{v}_n$.     It is easy to see that for each evaluation $\theta$,   the   Bellman recursive formula can be written   as follows:  $$ v_{\theta} (z)= \sup_{z'\in F(z)} \left( \theta_1 r(z')+ (1-\theta_1) \; v_{\theta^+}(z')\right),$$
\noindent where if $\theta_1<1$,      the ``shifted" evaluation $\theta^+$ is defined as $(\frac{\theta_{t+1}}{1- \theta_1})_{t \geq 1}$. 

\vspace{0,5cm}
 \begin{lem} \label{lem1} For all evaluation  $\theta$ in $\Theta$ and state $z$ in $Z$,
$$  |v_{\theta}(z) -\sup_{z' \in F(z)} v_{\theta} (z')|\leq \theta_1+ \sum_{t\geq 2} | \theta_t- \theta_{t-1}|.$$\end{lem}    
    
    \noindent{\bf Proof:} Consider any $z_1\in F(z)$,  and for   $\varepsilon>0$ a play $s=(z_2,z_3,...,)$ in $S(z_1)$ such that $\gamma_\theta(s)\geq v_\theta(z_1)-\varepsilon$. We have:
    \begin{eqnarray*}
    v_\theta(z) & \geq & \theta_1 r(z_1) + \sum_{t=2}^\infty \theta_t r(z_t)\\
     & \geq & \theta_1 r(z_1) + \sum_{t=2}^\infty \theta_{t-1} r(z_t) + \sum_{t=2}^\infty (\theta_{t}-\theta_{t-1}) r(z_t)\\
     & \geq &  v_\theta(z_1) -\varepsilon -  \sum_{t=2}^\infty | \theta_{t}-\theta_{t-1}|.
     \end{eqnarray*}
     Conversely, choose $s=(z_1,z_2,...)$ in $S(z_0)$ such that $\gamma_\theta(s)\geq v_\theta(z)-\varepsilon$.
       \begin{eqnarray*}
    v_\theta(z) & \leq & \varepsilon +  \theta_1 r(z_1) + \sum_{t=2}^\infty \theta_{t-1} r(z_t) + \sum_{t=2}^\infty (\theta_{t}-\theta_{t-1}) r(z_t)\\ 
     & \leq &   \varepsilon+ \theta_1 + v_\theta(z_1) +  \sum_{t=2}^\infty | \theta_{t}-\theta_{t-1}|.  
     \end{eqnarray*}
$\;$ \hfill $\Box$

  \begin{defi}   The total variation of an evaluation $\theta=(\theta_t)_{t\geq 1}$  is $$TV(\theta)=\sum_{t=1}^\infty |\theta_{t+1}-\theta_t|.$$ \end{defi}
    
 \noindent We have $\sup_t \theta_t \leq TV(\theta) \leq 2$. 
In the case of a Ces\`{a}ro evaluation  $\theta=(1/n,...,1/n,$ $0,0,...)$, we have $TV(\theta)=1/n$. For a discounted evaluation $\theta=(\lambda {(1-\lambda)}^{t-1})_{t \geq 1}$, we have $TV(\theta)=\lambda$.  A small  $TV(\theta)$ corresponds to  a patient evaluation, and sometimes we will refer to $TV(\theta)$ as the {\it  impatience} of $\theta$. We will consider here limits  when $TV(\theta)$ goes to zero, generalizing the cases where $n \longrightarrow  \infty$ or $\lambda \longrightarrow  0$.    Notice that if an evaluation $\theta$ is non increasing, i.e. satisfies $\theta_{t+1}\leq \theta_t$ for all $t$, we have that $TV(\theta)=\theta_1$. In the case of  a sequence of non increasing evaluations $(\theta^k)_k$,  the condition $TV(\theta^k) \xrightarrow[k \to \infty]{}0$ is equivalent to the condition $\sup_{t\geq 1} \theta_t^k  \xrightarrow[k \to \infty]{}   0$. We always have: 
$$(1-\theta_1) \sum_{t\geq 1} |\theta_t- \theta^+_t|\leq \theta_1 + TV(\theta),$$
\noindent so if $TV (\theta)$ is small,  the $L^1$-distance between $\theta$ and the shifted evaluation $\theta^+$ is also small. Notice also the following inequalities: for any given $T$, denote by $\overline{\theta}(T)$ the arithmetic mean  of $\theta_1$,..., $\theta_T$. We have for all $t=1,...,T$:
$$ |\theta_t -\overline{\theta}(T)|\leq \sum_{t'=1}^{T-1} |\theta_{t'}-\theta_{t'+1}|\leq TV(\theta).$$   
\noindent So if $TV(\theta)$ is small, then for all $T$ and $t\leq T$, the weight $\theta_t$ is close to  the average $\overline{\theta}(T)$. \\

   Given an evaluation $\theta$ and $m\geq 0$,   we write ${{v}}_{m, \theta}$ for the value function associated to the evaluation $\theta'= \sum_{t=1}^\infty  \theta_t \delta_{m+t}.$ The following  function will play a very important role in the sequel: 

\begin{defi} \label{defv^*}Define for all $z$ in $Z$, $$v^*(z)=\inf_{\theta \in \Theta} \; \sup_{m \geq 0} \;  {{v}}_{m, \theta}(z).$$ \end{defi}


%

\subsection{Main results}

We now state the main result of this paper. Recall that a metric space   is   totally bounded, or precompact, if for all $\varepsilon>0$ it can be covered by finitely many balls with radius $\varepsilon$.

        \begin{thm}  \label{thm1} Let $(\theta^k)_{k \geq 1}$ be a sequence of evaluations such that $TV(\theta^k) \xrightarrow[k \to \infty]{} 0.$ We have for all $z$ in $Z$: 
        $$v^*(z)= \inf_{k\geq 1} \sup_{m \geq 0} v_{m,\theta^k}(z).$$
    Moreover,    the sequence $(v_{\theta^k})_k$ uniformly converges if and only if the metric space $(\{v_{\theta^k}, k \geq 1\}, d_\infty)$   is totally bounded. And in  case of convergence, the limit value  is  $v^*$.\end{thm}
        
  This theorem generalizes theorem 3.10 in Renault, 2011, which was only dealing with Ces\`{a}ro evaluations\footnote{In this paper it is also proved that if the Ces\`{a}ro values $(v_{\bar{n}})_n$ uniformly converge then the limit can only be    $ \inf_{n\geq 1} \sup_{m \geq 0} v_{m,\bar{n}}$ which in this case is also equal to $\sup_{m \geq 0}\inf_{n\geq 1}  v_{m,\bar{n}}$.}. 
     In particular, there is a unique possible limit point for all sequences $(v_{\theta^k})_k$   such that $TV(\theta^k) \xrightarrow[k \to \infty]{} 0$, and consequently   any (uniform) limit  of such  sequence     is $v^*$.   Notice that this  is not true if we replace uniform  convergence by pointwise convergence: even for uncontrolled problems, it may happen that several limit points are possible.  As an immediate corollary of theorem \ref{thm1},  when $Z$ is finite the sequence  $(v_{\theta^k})_k$ is bounded and has a unique limit point, so converges to $v^*$. 
              
        \begin{cor}  \label{cor1}   Assume that $Z$ is endowed with a distance $d$ such that:  a) $(Z,d)$ is a   precompact metric space, and b)  the family $(v_\theta)_{\theta \in \Theta}$  is   uniformly equicontinuous. Then there is general uniform convergence of the value functions to $v^*$, i.e. $$\forall \varepsilon >0, \exists \alpha >0, \forall \theta \in \Theta \; s.t. \; TV(\theta)\leq \alpha, \; \|v_\theta-v^*\| \leq \varepsilon.$$ \end{cor}
 
 The proof of corollary \ref{cor1} from theorem \ref{thm1} follows from 1) Ascoli's theorem, and 2) the fact that the convergence of   $(v_{\theta^k})_k$   to  $v^*$ for each sequence of evaluations such that $TV(\theta^k)\xrightarrow[k \to \infty]{} 0$ is enough to have the general uniform convergence of the value functions to $v^*$.  \\

  \begin{cor} \label{cor2} 
   Assume that  $Z$ is endowed with a distance $d$ such that: a)  $(Z,d)$ is a    precompact metric space, b)   $r$ is   uniformly continuous,  and c) $F$   is   non expansive, i.e. $\forall z\in Z, \forall z' \in Z, \forall z_1 \in F(z), \exists z'_1 \in F(z') \; s.t. \; d(z_1,z'_1)\leq d(z,z').$ Then we have the same conclusions as corollary \ref{cor1}, there is general uniform convergence of the value functions to $v^*$, i.e. $$\forall \varepsilon >0, \exists \alpha >0, \forall \theta \in \Theta \; s.t. \; TV(\theta)\leq \alpha, \; \|v_\theta-v^*\| \leq \varepsilon.$$ 
 \end{cor}

\noindent{\bf Proof of corollary \ref{cor2}.} One can proceed as in the proof of corollary 3.9 in Renault, 2011. Given two   states $z$ and $z'$, one can construct inductively from each play $s=(z_t)_{t \geq 1}$ at $z$ a play $s=(z'_t)_{t \geq 1}$ at $z'$ such that $d(z_t,z'_t)\leq d(z,z')$ for all $t$.  Regarding payoffs, we introduce the modulus of continuity $\hat{\varepsilon}$ of $r$ by:  
\centerline{$\hat{\varepsilon}(\alpha)= \sup_{ z, z'\rm s.t. \it d(z,z')\leq \alpha} |r(z)-r(z')|$ for each    $\alpha\geq 0$. } So $|r(z)-r(z')|\leq \hat{\varepsilon}(d(z,z'))$ for each pair of states $z$, $z'$, and  $\hat{\varepsilon}$ is continuous at 0. Using the previous construction, we obtain that for $z$ and $z'$ in $Z$,  for all  $k\geq 1$,  $ |v_{\theta^k}(z)-v_{\theta^k}(z')|\leq \hat{\varepsilon}(d(z,z'))$.  In particular, the family $(v_{\theta^k})_{k\geq 1}$ is uniformly continuous, and  corollary  \ref{cor1} gives the result. \hfill $\Box$

\vspace{0,3cm}

 A  completely different  proof of corollary \ref{cor2}, with another expression for the limit value $v^*$, can be found in theorem 3.9 of Renault Venel 2012.

 \vspace{1cm}
 
\subsection {Extension to stochastic transitions}
 We generalize here theorem \ref{thm1} to the case of stochastic transitions. We will only consider transitions with finite support, and given a set $X$ we denote by $\Delta_f(X)$ the set of probabilities  with finite support over $X$.  We consider now   stochastic dynamic programming problems of the following form. There is an arbitrary  non empty set of states $X$,  a transition given by a multi-valued mapping  $F:X  \rightrightarrows \Delta_f(X)$ with non empty values, and a payoff (or reward) function $r:X\rightarrow [0,1]$.  The interpretation  is that given an initial state  $x_0$ in $X$, a decision-maker    has to  choose a probability with finite support $u_1$ in $F(x_0)$, then $x_1$ is selected according to $u_1$ and there is a payoff $r(x_1)$. Then the player has to  select $u_2$ in $F(x_1)$, $x_2$ is selected according to $u_1$ and the player receives the payoff $r(x_2)$, etc... 
 
 Following Maitra and Sudderth (1996), we say that  $\Gamma=(X, F, r)$ is a Gambling House.  We assimilate an  element $x$ in $X$ with its  Dirac measure $\delta_x$  in $\Delta(X)$, we write $Z=\Delta_f(X)$ and an element in $Z$ is written $u=\sum_{x \in X} u(x) \delta_x$. In case  the values of $F$ only consist of  Dirac measures on $X$, we are in the previous case of a dynamic programming problem.
 
 We linearly extend $r$ and $F$   to $\Delta_f(X)$ by defining for each $u$ in $Z$,   the payoff $r(u)=\sum_{x \in X} r(x) u(x)$ and the transition $F(u)=\{\sum_{x \in X} u(x)f(x), s.t. \ f:X \rightarrow Z \text{ and } f(x)\in F(x) \  \forall x \in X\}$. A  play  at $x_0$ is a sequence $\sigma=(u_1,...,u_t,...)\in Z^\infty$ such that  $u_1\in F(x_0)$ and $u_{t+1} \in F(u_{t})$ for each $t\geq 1$, and we denote by $\Sigma(x_0)$ the set of plays at $x_0$. Given  an evaluation  $\theta$,  the $\theta$-payoff of a  play $\sigma=(u_1,...,u_t,...)$ is defined as: $\gamma_{\theta}(\sigma)=\sum_{t\geq 1} \theta_t r(u_t)$, and the $\theta$-value at $x_0$ is:
 $$v_{\theta}(x_0)=\sup_{\sigma \in \Sigma(x_0)}  \gamma_{\theta}(\sigma).$$
  \noindent $v_\theta$ is by definition a mapping from $X$ to $[0,1]$, and we define as before, for all $x$ in $X$:
  $$v^*(x)=\inf_{\theta \in \Theta} \; \sup_{m \geq 0} \;  {{v}}_{m, \theta}(x).$$
 Theorem 1 easily  extends to this context.

    \begin{thm}  \label{thm2} Let $(\theta^k)_{k \geq 1}$ be a sequence of evaluations with vanishing total variation, i.e. such that $TV(\theta^k) \xrightarrow[k \to \infty]{} 0.$ We have:
        $$\forall x \in X, \; \; v^*(x)= \inf_{k\geq 1} \sup_{m \geq 0} v_{m,\theta^k}(x).$$
    Moreover,    the sequence $(v_{\theta^k})_k$ uniformly converges if and only if the metric space $(\{v_{\theta^k}, k \geq 1\}, d_\infty)$   is totally bounded. And in case of convergence, the limit value  is  $v^*$.\end{thm}

 \noindent{\bf Proof.} Consider the deterministic dynamic programming problem $\Gamma=(Z,F,r)$. For any evaluation $\theta$, the associated $\theta$-value function $\tilde{v}_\theta: Z \longrightarrow [0,1]$ is the affine extension of ${v}_\theta: X \longrightarrow [0,1]$. We put,  as in definition \ref{defv^*}, for all $z$ in $Z$:
 $$\tilde{v}^*(z)=\inf_{\theta \in \Theta} \; \sup_{m \geq 0} \;  {\tilde{v}}_{m, \theta}(z).$$
\noindent Notice that as an ``infsup" of affine functions,   there is no reason a priori  for $\tilde{v}^*$ to be affine. However, the restriction of $\tilde{v}^*$ to $X$ is $v^*$. 

Consider now a sequence $(\theta^k)_{k \geq 1}$ of evaluations with vanishing total variation. Applying theorem \ref{thm1} to $\Gamma$, we first obtain that for all $x$ in $X$: $$v^*(x)= \inf_{k\geq 1} \sup_{m \geq 0} v_{m,\theta^k}(x).$$
\noindent Moreover, given two evaluations $\theta$ and $\theta'$, we have (using  the same notation $d_\infty$ for the distances on $[0,1]^X$ and on $[0,1]^Z$):
 \begin{eqnarray*}
d_{\infty} (\tilde{v}_{\theta}, \tilde{v}_{\theta'}) & = & \sup_{z \in Z} |\tilde{v}_\theta(z)- \tilde{v}_{\theta'}(z)|,\\
\;  & = &   \sup_{z \in Z} |\int_{p \in X} {v}_{\theta}(p)-  {v}_{\theta'}(p) du(p)|,  \\
\;   & = & d_{\infty} ( {v}_{\theta},  {v}_{\theta'}). 
\end{eqnarray*}
\noindent Consequently, $(\{v_{\theta^k}, k \geq 1\}, d_\infty)$   is totally bounded if and only if $(\{v_{\theta'^k}, k \geq 1\}, d_\infty)$ is, and this completes the proof. \hfill $\Box$.

 \section{Examples}
 
The first  very simple example   shows that, even when the set of states is finite,  it  is not possible to obtain the   conclusions of theorem  \ref{thm1} or corollaries \ref{cor1} and \ref{cor2} with sequences of evaluations satisfying the  weaker convergence condition:  $\sup_{t\geq 1} \theta_t^k \longrightarrow_{k \to \infty} 0$.  

\begin{exa} \rm Consider   the following  dynamic  programming problem with 2 states: $Z=\{z_0,z_1\}$, $F(z_0)=\{z_1\}$, $F(z_1)=\{z_0\}$, with payoffs $r(z_0)=0$ and $r(z_1)=1$. We have a deterministic Markov chain, so that any play alternates forever between $z_0$ and $z_1$. Define for each $k$ the evaluations $\theta^k=\frac{1}{k} \sum_{t=1}^k \delta_{2t-1}$ and $\theta'^k=\frac{1}{k} \sum_{t=1}^k \delta_{2t}$. We have $v_{\theta^k}(z_0)=v_{\theta'^k}(z_1)=1$, and $v_{\theta^k}(z_1)=v_{\theta'^k}(z_0)=0$ for all $k$. Define now $\nu^k$ as $\theta^k$ when $k$ is even, and $\theta'^k$ when $k$ is odd. The evaluation $\nu^k$ satisfies $\sup_t \nu^k_t=\frac{1}{k}\longrightarrow_{k \to \infty} 0$, however $(v_{\nu^k}(z_0))_k$ and $(v_{\nu^k}(z_1))_k$ do not converge. \hfill $\Box$
\end{exa}

\vspace{0,5cm}

Lehrer and Sorin (1992) proved that the uniform convergence of the Ces\`{a}ro values $(v_{\overline{n}})_{n \geq 1}$ was equivalent to the uniform convergence of the discounted values $(v_{\lambda})_{\lambda \in (0,1]}$. The following example shows that this property does not extend to general evaluations:  given 2 sequences of TV- vanishing evaluations $(\theta^k)_{k \geq 1}$ and $(\theta'^k)_{k \geq 1}$, the uniform convergence of   $(v_{\theta^k})_k$ and $(v_{\theta'^k})_k$ are not equivalent. 

\begin{exa} \label{exa2} \rm   In this example, $({v}_{\bar n})_n$ will pointwise  converges to the constant 1/2 whereas for a particular sequence of evaluations ($\theta^k)_k$ with total variation going to zero, we will have   $(v_{\theta^k})_k(z)= 1$ for all $k$ and $z$. 

We construct  a dynamic programming  problem defined via  a rooted tree $T$ without terminal nodes (as in Sorin Monderer 1992 or Lehrer Monderer  1994). $T$ has countably many nodes, and the payoff attached to each node is either 0 or 1.

We first construct  a tree $T_1$, with   countably many nodes and   root $z_0$. Each node has an outgoing degree one, except the root which has countably many potential successors $z_1$, $z_2$,..., $z_n$... On the $n^{th}$ branch starting from $z_n$, each node has a unique successor and the payoffs starting from $z_n$ are successively 0 for $n$ stages, then 1 for $n$ stages, then 0 until the end of the play.

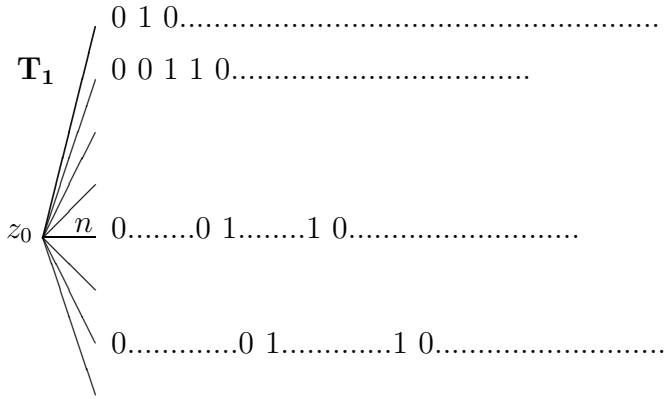
\unitlength 0,7mm

 \begin{center}
 \begin{picture}(200,100)
 \put(0,80){${\bf T_1}$}
 \put(-2,50){$z_0$}

 \put(5,50){\line(1,4){10}}
  
 \put(5,50){\line(1,3){10}}
  \put(5,50){\line(1,2){10}}
   \put(5,50){\line(1,1){10}}
  
 \put(5,50){\line(1,0){10}}
  
 \put(5,50){\line(1,-1){10}}
  \put(5,50){\line(1,-2){10}}
  \put(5,50){\line(1,-3){10}}
   \put(5,50){\line(1,4){10}}

  \put(11,51){$n$}
    
     \put(18,90){$0$ $1$ $0$........................................................}
     
          \put(18,80){$0$ $0$ $1$ $1$ $0$...................................}
 
    \put(18,50){$0$........$0$ $1$........$1$ $0$...........................}
 
   \put(18,28){$0$.............$0$ $1$.............$1$ $0$...........................}

\end{picture}
  \end{center}

We now define $T$ inductively from $T_1$. $T_2$ is obtained from $T_1$ by attaching the tree $T_1$ to each node of $T_1\backslash\{z_0\}$. This means that  for  each node $z$ of $T_1\backslash\{z_0\}$ we add a copy of the tree $T_1$ where $z$ plays the role of the root of $T_1$.  And for each $l$, the tree $T_l$  is obtained   by attaching the tree $T_1$ to each node of $T_{l-1}\backslash T_{l-2}$. Finally, $T$ is defined as the union $\bigcup_{l\geq 1} T_l$.

Starting from $z_0$, any sequence of $n$ consecutive payoffs of 1  has to be preceeded by $n$ consecutive payoffs of 0, so $\bar{v}_n(z_0)\leq 1/2$ for each $n\geq 1$, and for each node $z$ and even integer $n$ it is possible to get exactly $n/2$ payoffs of 0 followed by $n/2$ payoffs of 1. Consequently one can deduce that $(v_n(z))_n$ converges to 1/2 for each state $z$. But $\sup_{z\in Z} v_n(z)=1$ for each $n$, and the convergence is not uniform.

Consider now for any $k$, the evaluation $\theta^k=(0,...0,\frac{1}{K},..., \frac{1}{K},0,...)$ $=\frac{1}{K}\sum_{t=1}^K \delta_{t+K}.$ We have $v_{\theta^k}(z)=1$ for all $k$ and $z$, so $(v_{\theta^k})_k$ uniformly converges to $v^*=1$. \hfill $\Box$
  \end{exa}

\begin{exa}\rm The condition $(\{v_{\theta}, \theta \in \Theta\}, d_\infty)$    totally bounded is satisfied with the hypotheses of corollary \ref{cor1} or  corollary \ref{cor2}, and is sufficient to obtain the general uniform convergence of the value functions. This condition turns out to be stronger than having $(\{v_{\theta^k}, k \geq 1\}, d_\infty)$ totally bounded for every sequence of evaluations with vanishing TV. 

In the following example, there is no control and the state space  $Z$ is  the set of all integers, with transition given by the shift: $F(z)=\{z+1\}$. The payoffs are given by $r(0)=1$ and $r(z)=0$ for all $z\neq 0$.

 For all evaluations $\theta=(\theta_t)_{t \geq 1}$, we have $\sup_{z \in Z} v_\theta(z)=\sup_t \theta_t$, so we have general uniform convergence of the value functions to $v^*=0$.

For all positive $t$, we can consider the evaluation given by the Dirac measure on $t$. We have $v_{\delta_t}(-t)=1$, and $v_{\delta_t}(z)=0$ if $z\neq -t$. The set $\{v_{\delta_t}, t \geq 1\}$ is not totally bounded. \hfill $\Box$
\end{exa}

\section{Proof of theorem \ref{thm1}}

We start with a few notations and definitions. 
We  define inductively a sequence of correspondences $(F^n)_n$ from $Z$ to $Z$, by   $F^0(z)=\{z\}$ for every state $z$, and $\forall n\geq 0$, $F^{n+1}=F^n\circ F$ (the composition being  defined by $G\circ H(z)=\{z" \in Z, \exists z'\in H(z), z"\in G(z')\}$). $F^n(z)$  represents the set of states that the decision maker can reach in $n$ stages from the initial state  $z$. We also define  for every  state $z$,  $G^m(z)= \bigcup_{n=0}^m  {F}^n(z)$ and $G^{\infty}(z)=  {\bigcup_{n=0}^{\infty} F^n(z)}$. The set $G^{\infty}(z)$ is  the set of states that the decision maker,  starting from $z$,  can reach in a finite number of stages. 

For all $\theta$ in $\Theta$, $m\geq 0$ and  initial state  $z$, we clearly have: 
$$v_{m, \theta}(z)=\sup_{z'\in F^m(z)} v_\theta(z')=\sup_{s \in S(z)} \sum_{t=1}^\infty \theta_t r(z_{m+t}).$$



%

\noindent In the sequel, we fix  a sequence of evaluations $(\theta^k)_{k \geq 1}$ such that $TV(\theta^k) \xrightarrow[k \to \infty]{} 0.$

 \begin{lem} \label{lem2} For all $m_0 \geq 0$ and $z$ in $Z$,
 $$\liminf_k  \sup_{m \leq m_0} v_{m, \theta^k}(z)=\liminf_k v_{\theta^k}(z).$$
  \end{lem}
\noindent{\bf Proof:} For each $k$, we have $v_{\theta^k}(z)\geq v_{1, \theta^k}(z) - \theta^k_1 - TV(\theta^k)$ by lemma \ref{lem1}, so  $v_{\theta^k}(z)\geq v_{1, \theta^k}(z) -  2 TV(\theta^k)$. Iterating, we obtain that:

 $v_{\theta^k}(z)\geq\sup_{m \leq m_0} v_{m, \theta^k}(z) - 2 m_0 TV(\theta^{k}).$ \hfill $\Box$

\vspace{0,5cm}

 A key result is the following proposition, which is true for all evaluations $\theta$.    \begin{pro} \label{pro1} For all evaluations $\theta$ in $\Theta$ and initial state $z$ in $Z$,
    $$\sup_{z'\in G^\infty(z)} v_\theta(z')\geq \limsup_k {v_{\theta^k}}(z).$$ \end{pro}
    
    \noindent{\bf Proof of proposition  \ref{pro1}} $z$ and $\theta$ being fixed, put $\beta=\sup_{z'\in G^{\infty}(z)} v_\theta(z')$. Fix $\varepsilon \in (0,1]$, there exists $T_0$ such that $\sum_{t=T_0+1}^\infty \theta_t \leq \varepsilon$, and fix  $T_1\geq T_0/\varepsilon$.      

For any play $s=(z_1,...,z_t,...)$ in $S(z)$, we have by definition of $\beta$ that for all $T$, $\sum_{t =T+1}^\infty \theta_{t-T} r(z_t)\leq \beta$. Let $m$ be a non negative integer,   we  define: $$A_m=\sum_{T=m T_1}^{(m+1)T_1-1} \sum_{t=T+1}^\infty \theta_{t-T} r(z_t) \leq  T_1  \beta.$$
\begin{eqnarray*}
A_m &=& \sum_{t=mT_1+1}^\infty r(z_t)  \sum_{T=m T_1}^{\min\{(m+1)T_1-1,t-1\}} \theta_{t-T}, \\
 &\geq & \sum_{t=mT_1+1}^{(m+1)T_1} r(z_t) (\theta_1+...+\theta_{t-mT_1}), \\
 & \geq & (1-\varepsilon) \sum_{t=T_0+mT_1}^{(m+1)T_1} r(z_t),\\
  & \geq & (1-\varepsilon) \left(\sum_{t=1+mT_1}^{(m+1)T_1} r(z_t)- (T_0-1)\right).
  \end{eqnarray*}
We obtain:  $$ T_1 \beta \geq  (1-\varepsilon) T_1 \overline{\gamma}_{m T_1, T_1}(s) - (1-\varepsilon) (T_0-1),$$
and \begin{equation} \label{eq1} \overline{\gamma}_{mT_1, T_1}(s)\leq \frac{\beta}{1-\varepsilon} +  \varepsilon.\end{equation}

We now consider $\gamma_{\theta^k}(s)$ for $k$ large.  We compute $\sum_{t=1}^\infty \theta^k_t r(z_t)$ by dividing the stages into blocks of length $T_1$. For each $m\geq 0$, let $ \overline{\theta^k}(m)$ be the Ces\`{a}ro-average of $\theta^k_t$, where $t$ ranges from $mT_1+1$ to $(m+1)T_1$. Notice that for all such $t$, we have $|\theta^k_t- \overline{\theta^k}(m)|$ $\leq$ $\sum_{t'=mT_1+1}^{(m+1)T_1-1} |\theta_{t'}-\theta_{t'+1}|$. We have:
\begin{eqnarray*}
\sum_{t =mT_1+1}^{(m+1)T_1} \theta^k_t r(z_t) & \leq & \sum_{t =mT_1+1}^{(m+1)T_1}  \overline{\theta^k}(m)   r(z_t) + \sum_{t =mT_1+1}^{(m+1)T_1} |\theta^k_t - \overline{\theta^k}(m)|  r(z_t) ,\\
 & \leq  &\overline{\theta^k}(m)  \overline{\gamma}_{mT_1, T_1}(s) + T_1 \sum_{t =mT_1+1}^{(m+1)T_1-1} |\theta^k_t-\theta^k_{t+1}|,\\
  & \leq & \overline{\theta^k}(m)  \left( \frac{\beta}{1-\varepsilon} +  \varepsilon \right) + T_1 \sum_{t =mT_1+1}^{(m+1)T_1-1} |\theta^k_t-\theta^k_{t+1}|.
  \end{eqnarray*}
  \noindent where the last inequality follows from equation (\ref{eq1}). Summing up over $m$, we obtain:
  $$\gamma_{\theta^k}(s)\leq  \frac{\beta}{1-\varepsilon} +  \varepsilon + T_1 \; TV(\theta^k).$$
Consequently,  $ \limsup_k {v_{\theta^k}}(z)\leq  \frac{\beta}{1-\varepsilon} +  \varepsilon$, and this is true for all $\varepsilon$. 
  \hfill $\Box$

\vspace{1cm}

\begin{cor} \label{cor3} $$\inf_{\theta \in \Theta}\;  \sup_{m \geq 0} v_{m, \theta} = \inf_{k \geq 1} \sup_{m\geq 0} v_{m, \theta^k} .$$ \end{cor}

\noindent{\bf Proof:} Consider an initial state $z$, and write $\alpha= \inf_{k} \sup_{m} v_{m, \theta^k}(z)$.  It is clear that $\alpha  \geq$ $\inf_{\theta \in \Theta} \sup_{m \geq 0} v_{m, \theta}(z)$. Now for each $k\geq 1$ there exists $m(k)$, such that $v_{m(k), \theta^k}(z)\geq 
\alpha  -1/k$, and we define  the evaluation  $\theta'^{k}=\sum_{t=m(k)+1}^{\infty} \theta^k_{t-m(k)} \delta_t$. We have  $TV (\theta'^{k})= TV (\theta^k) \xrightarrow[k \to \infty]{} 0$, so by proposition  \ref{pro1} we obtain that for all evaluations $\theta$, $\sup_{z'\in G^\infty(z)} v_\theta(z')\geq \limsup_k {v_{\theta^{m(k),k}}}(z)$ $\geq  \alpha .$ \hfill $\Box$ 

\vspace{0,5cm}

From lemma \ref{lem2} and proposition \ref{pro1}, one can easily deduce the following corollary.

\begin{cor} \label{cor4} For all $m_0 \geq 0$ and $z$ in $Z$,
$$\inf_{k\geq 1}  \;  \sup_{m \leq m_0} v_{m, \theta^k}(z) \leq  \liminf_{k}   v_{\theta^k}(z) \leq \limsup_{k}  v_{\theta^k}(z) \leq \inf_{k \geq 1} \sup_{m\geq 0} v_{m, \theta^k}(z) .$$ \end{cor}

And we can now conclude the proof the theorem \ref{thm1}, proceeding as in the proof of theorem 3.10 in Renault, 2011.\\

\noindent{\bf End of the proof of theorem \ref{thm1}}.   Define $d(z,z')=\sup_{k\geq 1} |v_{\theta^k}(z)-v_{\theta^k}(z')|$ for all states $z$ and $z'$. The space $(Z,d)$ is now a pseudometric space (may not be Hausdorff). By assumption, there exists a finite set of indices $I$ such that for all $k\geq 1$, there exists $i$ in $I$ satisfying $d_{\infty}(v_\theta^k,v_ i)\leq \varepsilon$. Consider now the set $\{(v_i(z))_{i \in I}, z \in Z\}$, it is a subset of the compact metric space $[0,1]^I$ with the uniform distance, so it is itself precompact and we obtain the existence  of a finite subset $C$ of states in $Z$ such that:$$\forall z \in Z, \exists c \in C, \forall i \in I, |v_i(z)-v_i(c)|\leq \varepsilon.$$
We have obtained that for each $\varepsilon>0$, there exists a finite subset $C$ of $Z$ such that for every $z$ in $Z$, there is $c\in C$ with $d(z,c)\leq \varepsilon$. The pseudometric space ($Z,d)$ is itself precompact. Equivalently, any sequence in $Z$ admits a Cauchy subsequence  for $d$. Notice that all value functions $v_{\theta^k}$ are clearly 1-Lipschitz for $d$.\\

Fix $z$ in $Z$, and consider now the sequence of sets $(G^m(z))_{m\geq 0}$. For all $m$, $G^m(z)\subset G^{m+1}(z)$ so using the precompacity of $(Z,d)$ it is not difficult to show (see, e.g. step 2 in the proof of theorem 3.7 in Renault, 2011) that $(G^m(z))_{m\geq 0}$ converges to $G^\infty(z)$, in the sense that:
\begin{equation}\label{eq2} \forall \varepsilon>0, \exists m\geq 0, \forall z' \in G^\infty(z), \exists z''\in G^m(z), \; d(z',z'')\leq \varepsilon.\end{equation}

We now use corollary \ref{cor4} to conclude. We have for all $m$ : $$\inf_{k\geq 1}  \;  \sup_{z' \in G^{m}(z)} v_{\theta^k}(z') \leq  \liminf_{k}   v_{\theta^k}(z) \leq \limsup_{k}  v_{\theta^k}(z) \leq \inf_{k \geq 1} \sup_{z'\in G^\infty(z)} v_{ \theta^k}(z') .$$

\noindent Fix finally $\varepsilon>0$, and consider $k\geq 1$ and $m\geq 0$ given by equation (\ref{eq2}). Let $z'$   in $G^\infty(z)$ be such that $v_{\theta^k}(z')\geq  \sup_{z'\in G^\infty(z)} v_{ \theta^k}(z')-\varepsilon.$ Let $z''$  in $G^m(z)$ be such that $d(z',z'')\leq \varepsilon$. Since $v_{\theta^k}$ is  1-Lipschitz for $d$, we obtain $v_{\theta^k}(z'')\geq  \sup_{z'\in G^\infty(z)} v_{ \theta^k}(z')-2 \varepsilon.$ Consequently,  $ \sup_{z'\in G^m(z)} v_{ \theta^k}(z')\geq \sup_{z'\in G^\infty(z)} v_{ \theta^k}(z')-2 \varepsilon$ for all $k$, so 
$$\inf_{k\geq 1}  \;  \sup_{z' \in G^{m}(z)} v_{\theta^k}(z') \geq    \inf_{k \geq 1} \sup_{z'\in G^\infty(z)} v_{ \theta^k}(z') - 2 \varepsilon.$$
We obtain $ \liminf_{k \geq 1}   v_{\theta^k}(z) \geq \limsup_{k \geq 1}  v_{\theta^k}(z)-2 \varepsilon$, and so $(v_{\theta^k}(z))_k$ converges. Since $(Z,d)$ is precompact and all $v_{\theta^k}$ are 1-Lipschitz, the convergence is uniform.  $\Box$

\section{An open question}

We know since Lehrer and Sorin (1992)   that the uniform convergence of the Ces\`{a}ro values $(v_{\overline{n}})_{n \geq 1}$ is  equivalent to the uniform convergence of the discounted values $(v_{\lambda})_{\lambda \in (0,1]}$. Example \ref{exa2} shows that is possible to have no uniform convergence of the Ces\`{a}ro values  (or equivalently of the discounted values) but uniform convergence for a particular sequence of evaluations with vanishing TV. Could it be the case that the Ces\`{a}ro values and the discounted values have  the following ``universal" property ?  

Assuming  uniform convergence of the Ces\`{a}ro values, do we have general uniform convergence of the value functions, i.e. is it true that  $(v_{\theta^k})_k$ uniformly converges for every  sequence of evaluations $(\theta^k)_{k \geq 1}$ such that  $TV(\theta^k) \xrightarrow[k \to \infty]{} 0$ ? 

\vspace{0,5cm}

The above property is true in case of an uncontrolled problem (zero-player), i.e. when the transition $F$ is single-valued.

\begin{pro} For an uncontrolled problem, the uniform convergence of the   Ces\`{a}ro values implies the general uniform convergence of the value functions:  $$\forall \varepsilon >0, \exists \alpha >0, \forall \theta \in \Theta \; s.t. \; TV(\theta)\leq \alpha, \; \|v_\theta-v^*\| \leq \varepsilon.$$ \end{pro}

\noindent {\bf Proof:} Fix $\varepsilon>0$. By assumption there exists $N$ such that for all states $z$ in $Z$, $|\bar{v}_N(z)-v^*(z)|\leq \varepsilon.$ Consider an arbitrary evaluation $\theta$ and an initial state  $z_0$.  For each positive $t$  we denote by $z_t$ the state reached from $z_0$ in $t$ stages, we have $v_\theta(z_0) =\sum_{t=1}^\infty \theta_t r(z_t)$ and $v^*(z_0)=v^*(z_t)$ for all $t$.

Divide the set of stages into consecutive blocks of length $N$: $B^0=\{1,..., N\}$,..., $B^m=\{m N+1,..., (m+1) N\}$,... Denote by $\bar{\theta}(m)$ the mean of $\theta$ over $B^m$, we have $\sum_{m=0}^\infty N \bar{\theta}(m)=1$. We  also write $\bar{r}(m)$ for the mean $\frac{1}{N}\sum_{t \in B^m} r(z_t)$. We have $\bar{r}(m)= \bar{v}_{N}(z_{m N})$, so $|\bar{r}(m)-v^*(z_0)|\leq \varepsilon$ for all $m$. 

Computing payoffs by blocks, we have
\begin{eqnarray*}
 v_\theta(z_0) &=& \sum_{m=0}^\infty \sum_{t \in B^m} \theta_t r(z_t),\\
  & = &  \sum_{m=0}^\infty \sum_{t \in B^m}( \theta_t- \bar{\theta}(m)) r(z_t) + \sum_{m=0}^\infty  N \bar{\theta}(m) \bar{r}(m).
  \end{eqnarray*}
 \noindent So we obtain:
 $$  v_\theta(z_0)- v^*(z_0)= \sum_{m=0}^\infty \sum_{t \in B^m}( \theta_t- \bar{\theta}(m)) r(z_t) + \sum_{m=0}^\infty  N \bar{\theta}(m) ( \bar{r}(m)- v^*(z_0)),$$
 \noindent and 
$$  |v_\theta(z_0)- v^*(z_0)| \leq \sum_{m=0}^\infty  N \sum_{t \in B^m} |\theta_{t+1}-\theta_t| + \varepsilon\leq N \; TV(\theta) + \varepsilon. $$
If $TV(\theta)\leq \frac{\varepsilon}{N}$, we get $  |v_\theta(z_0)- v^*(z_0)| \leq 2 \varepsilon$, hence the result. \hfill $\Box$

\section{References}

   








\noindent Blackwell, D. (1962): Discrete dynamic programming. The Annals of  Mathematical Statistics, 33, 719-726. \\

\noindent  Lehrer, E. and D. Monderer    (1994): Discounting versus Averaging in  Dynamic Programming. Games and Economic Behavior, 6, 97-113.\\

\noindent Lehrer, E. and S. Sorin (1992): A uniform Tauberian Theorem in Dynamic Programming. Mathematics of Operations Research, 17, 303-307.\\


%

\noindent Lippman, S. (1969): Criterion Equivalence in Discrete Dynamic Programming. Operations Research 17, 920-923.  \\



\noindent Maitra, A.P. and Sudderth, W.D. (1996): Discrete gambling and stochastic games,  Springer Verlag\\



\noindent Mertens, J-F. and A. Neyman (1981):  Stochastic games. International Journal of Game Theory, 1, 39-64.\\



\noindent Monderer, D. and S. Sorin (1993): Asymptotic properties in Dynamic Programming. International Journal of Game Theory, 22, 1-11.\\


 
%

\noindent Renault, J. (2011): Uniform value in Dynamic Programming.    Journal of the European Mathematical Society,  vol. 13, p.309-330. \\

\noindent Renault, J. and X. Venel (2012): A distance for probability spaces, and  long-term   values
  in Markov Decision Processes and Repeated Games.  preprint hal-00674998

\end{document}